\documentstyle[12pt,amssymb]{article}
\author{Mohammad~Obiedat}
\title{ A note on the localization of $J$-groups }

\begin{document}
\begin{center}
 {\large \textbf{A note on the localization of $J$-groups}\vspace{0.2cm}\\
   Mohammad~Obiedat\\
(Received June 4, 1998)\\
(Revised September 14, 1998)}
\end{center}
\begin{abstract}
 Let  $\widetilde{JO} (X) = \widetilde{KO}(X) / TO(X)  $ be the $J$-group
of a connected finite $CW$ complex $X.$ Using Atiyah-Tall
\cite{ati-tal}, we obtain two computable formulae of $TO(X)_{(p)},$ the
localization of $TO(X)$ at a prime
$p.$ Then we show how to use those two formulae of $TO(X)_{(p)}$ to find the
$J$-orders of elements of $\widetilde{KO}(X),$ at least the $2$ and
$3$ primary factors of the canonical generators of 
 $\widetilde{JO} ({\Bbb C}P^m).$ Here ${\Bbb C}P^m$ is the complex
projective space.
\end{abstract} 
{\large \textbf{1. Introduction}}

Let  $\widetilde{JO} (X) = \widetilde{KO}(X) / TO(X)  $ be the $J$-group
of a connected finite  $CW$ complex $X,$ where 
$\widetilde{KO}(X)$ is the additive subgroup of the  $KO$-ring $KO(X)$
of elements of virtual dimension zero  and $TO(X)
=\{E-F \in \widetilde{KO}(X):S(E\oplus n)$ is fibre homotopy
 equivalent to $S(F\oplus n) \mbox{ for some } n \in {\Bbb N}\}.$
Let $\psi ^{k}$  be the Adams operations.
Then  Adams~\cite{ada} and Quillen~\cite{qui} showed that
$TO(X)=WO(X)=VO(X):$
\begin{eqnarray} \label{eq:wox} 
WO(X)=\bigcap_{f}\widetilde{KSO}(X)_f
\end{eqnarray} 
where the intersection runs over all functions $f:{\Bbb N} \rightarrow
 {\Bbb N}$ and $\widetilde{KSO}(X)_f=$\\
$<k^{f(k)}(\psi^k-1)(u):u \in \widetilde{KSO}(X)\mbox{ and } k\in {\Bbb N}>,$
and\\
--------------------
 \newpage
\begin{eqnarray} \label{eq:vox} 
VO(X)&= \{ & x \in \widetilde{KSO}(X):
\mbox{ there exists } u \in \widetilde{KSO}(X) \mbox{ such that }
 \nonumber \\ & & \theta_{k}(x)
 =\frac{\psi^{k}(1+u)}{1+u} \mbox{ in } 1+\widetilde{KSO}(X)
 \otimes {\Bbb{Q}}_{k} \mbox{ for 
all } k \in  {\Bbb{N}}\} 
\end{eqnarray}
where 
$\theta_{k}$ are the Bott exponential classes, and
 ${\Bbb{Q}}_{k}=\{ n/k^m: n,m \in {\Bbb{Z}} \}.$

 For a prime $p,$ let  $\widetilde{JO} (X)_{(p)}$ denote the localization of
 $\widetilde{JO} (X)$ at $p.$ Since  $\widetilde{JO} (X)$ is a finite abelian
group, $\widetilde{JO} (X)_{(p)}$ is isomorphic to the $p$-summand of
$\widetilde{JO} (X)$. Moreover, since the
localization is an exact functor on the category of finitely generated 
abelian groups,   
$\widetilde{JO} (X)_{(p)} \cong \widetilde{KO}(X)_{(p)} / TO(X)_{(p)}.$
Using Atiyah-Tall \cite{ati-tal} we obtain 
two computable formulae of $TO(X)_{(p)}.$ The significance of those 
two localized formulae of $TO(X)$ is shown to find the $J-$orders
of elements of $\widetilde{KO}({\Bbb C}P^m).$ 

In \S2 using the fact that $\widetilde{KSO}(X)$ is an orientable
$\gamma$-ring and the $p$-adic completion $\widetilde{KSO}(X)_p$ 
is an orientable $p$-adic $\gamma$-ring, we define
 a natural exponential map 
$\theta_{k}^{or}:KSO(2)(X)\rightarrow KSO(X)_p$ for each positive integer $k.$ 
If $k$ is odd, $\theta_{k}^{or}$ is the 
extension of $\theta_{k}^{or}:VectSO(2)(X)$ $\rightarrow KSO(X)$ defined
in Dieck \cite{die}. From the main theorem of 
\cite{ati-tal}, we obtain the  commutative diagram in Theorem~{$2.3.$}

Our main result  is the following two formulae of
$TO(X)_{(p)},$ which can be obtained directly from Theorem~{2.3.} 
$$ TO(X)_{(p)}=(\psi^{k_p}-1)(\widetilde{KSO}(X)_{(p)}). \hspace{6.4cm} 
\mbox{(Formula~I)}$$  
\begin{eqnarray}
TO(X)_{(p)}=&\{&x \in 
\widetilde{KSO}(X)_{(p)}:\theta_{k_p}^{or}(x)=\frac{\psi^{k_p}(1+u)}{1+u} 
\mbox{ in } 1+\widetilde{KSO}(X)_{p} \nonumber \\ & & \mbox{ for some } u 
\in \widetilde{KSO}(X)_{p}\}. \hspace{4.7cm} \mbox{(Formula~II}) \nonumber
\end{eqnarray} 
Formula~I (resp. Formula~II) of $TO(X)_{(p)}$ may be thought of as the 
localization of $WO(X)$ (resp. $VO(X)$) at $p$.

Let $y=r \xi_m({\Bbb C}) -2 $  where $\xi_m({\Bbb C})$
  is the complex Hopf line 
bundle over ${\Bbb C} P^m.$ In \S3 we apply Formulae~I and II of
$TO(X)_{(p)}$ to find $b_m(P_{m}(y;m_1,\ldots,m_t)),$ the $J$-order
of $ P_{m}(y;m_1,\ldots,m_t)=m_1y+m_2 y^2+\cdots+m_ty^t \in
\widetilde{KO}({\Bbb C}P^m).$ \"Onder \cite{ond1} has given the formula
$TO(X)_{(2)}=(\psi^3-1)(\widetilde{KO}(X)_{(2)})$ by using Dieck \cite{die}
 Ch.~11, and applied this formula  to give 
 computation of the 2-primary factor of $b_m(y).$ We obtain sharper results in giving a simple 
formula for the 2 and 3 primary factors of the $J$-orders of the canonical
generators of $\widetilde{JO}({\Bbb C}P^m).$ 
Finally in \S4 we show how  Formulae I and II of $TO(X)_{(p)}$ can be
used to compute the group $\widetilde{JO}(X)$ for 
 our illustrative example  $X={\Bbb C}P^4.$\vspace{1cm}\\
{\large \textbf{2. Two computable formulae of $TO(X)_{(p)}$}}

Let $G$ be a finitely generated abelian group. For a prime $p$ let 
$G_{(p)}=\{ g/m:g\in G \mbox{ and } m \in {\Bbb Z} \mbox{ with } (p,m)=1\}$
denote the localization of $G$ at $p,$ then $G_{(p)}$ is canonically 
isomorphic to ${\Bbb Z}_{(p)} \otimes G.$
 Also, let $G_p=\lim_{\stackrel{\longleftarrow}{n}} G/p^nG$
 denote the $p$-adic completion of $G.$ Then $G_p$ is canonically isomorphic
to ${\Bbb Z}_p \otimes G.$ For a rational number $q$, 
$\nu_p(q)$  denotes the exponent of $p$ in the prime factorization of 
$q.$

L{\scriptsize{EMMA}} 2.1.
(i) {\it Let $G$ be a finite abelian group. Then the
 following groups are canonically isomorphic: 
$$G_{(p)}\cong G_p \cong G(p)$$
where $G(p)=\{g\in G:g \mbox{ has order a power of } p \}.$ Consequently, if
$g\in G $ has order $m,$ then the order of $g/1$ in $G_{(p)}=$ the order of
$1\otimes g \mbox{ in } G_p$ is equal to $p^{\nu_p(m)}.$}\\
(ii) {\it If $G$ is a finitely generated abelian gruop, then $G_{(p)}$
 is canonically embedded in $G_p.$} 

P{\scriptsize{ROOF}}. The proof is obvious.

Now, our aim is to show how to apply the work of Atiyah-Tall \cite{ati-tal}
 to find $\widetilde{JO}(X)_{(p)}.$ Let $KSO(d)(X)$ be the group obtained by
 symmetrization of the semigroup $VectSO(d)(X)$ of all isomorphic classes 
of real vector bundles over $X$ with structural group $SO(dn)$ for 
$n=1,2,\ldots$ . $KSO(d)(X)$ is monomorphically embedded in $KO(X)$ as the
subgroup of classes $x$ such that $\omega_1(x)=0$ and $\mbox{dim}(x)=dn$ 
 for some $n\in {\Bbb{N}},$ i.e., $KSO(d)(X)=\{E-F \in KO(X): 
\mbox{ dim}(E-F)=dn \mbox{ and } E,F \mbox{ are orientable}\}.$

Let $\widetilde{KSO}(d)(X)=\{E-F\in KSO(d)(X):\mbox{ dim }
 E=\mbox{ dim } F\}.$
It is easy to see that $KSO(d)(X)=d{\Bbb Z}\oplus \widetilde{KSO}(d)(X)$ and
$\widetilde{KSO}(d)(X)=\widetilde{KSO}(1)(X)$ for each $d\geq 1.$ So, for
simplicity, we write $\widetilde{KSO}(X)$ instead of $\widetilde{KSO}(d)(X).$
 It is well known that $\widetilde{KSO}(X)$ is an orientable $\gamma-$ring
 and $\widetilde{KSO}(X)_p$ is an orientable $p$-adic $\gamma-$ring
(see \cite{ati-tal}, or \cite{die} Ch. 3).

Let $k$ be an odd integer and $J$  be a set of $k$th roots of unity $u\neq 1$
which contains from each pair $u,u^{-1}$ exactly one element. The operations
$\theta_{k}^{or}:VectSO(2)(X)$\\
 $\rightarrow KSO(X)$ are defined in \cite{die}
and given by \begin{equation} \label{eq:thor}
\theta_{k}^{or}(E)=k^m\prod_{u\in J}\lambda_{-u}(E)
(1-u)^{-2m}=\prod_{u\in J}\lambda_{-u}(E)(-u^{-1})^m \end{equation}
  where $2m=\mbox{dim}\,E.$ $\theta_{k}^{or}$ does not depend on the
choice of $J$ \cite{die}.

If $(k,p)=1$ then $\theta_{k}^{or}(E)$ is invertible in $KSO(X)_p.$ So 
$\theta_{k}^{or}$ can be extended to $KSO(2)(X)$ with values in  $KSO(X)_p.$
 Also, by using the fact 
that $\theta_{k}^{or}$ is a natural exponential map, it can be shown that
 $\theta_{k}^{or}:\widetilde{KSO}(X)\rightarrow 1+\widetilde{KSO}(X)_p$
 where $1+\widetilde{KSO}(X)_p$ is the multiplicative group of elements
$1+w$ with $w\in \widetilde{KSO}(X)_p.$ The operations 
$\rho_{k}^{or}:\widetilde{KSO}(X)_p\rightarrow 1+\widetilde{KSO}(X)_p$ are
given by \begin{equation} \label{eq:rhor}
\rho_{k}^{or}(x)=\prod_{u\in J}\gamma_{u/u-1}(x). \end{equation}

Now we shall show how to define $\theta_{2k}^{or}$ and $\rho_{2k}^{or}$ for
$k\geq 1.$ If $p\neq 2$ then $1/2 \in {\Bbb Z}_p,$ So we can difine
 $\theta_{2k}^{or}:\widetilde{KSO}(X)\rightarrow 1+\widetilde{KSO}(X)_p.$ by
\begin{equation} \label{eq:thoe}
\theta_{2k}^{or}(E-F)=(\prod_{u^{2k}-1=0 \atop u\neq 1}\lambda_{-u}(E)
(\prod_{{u^{2k}-1=0} \atop {u\neq 1}}\lambda_{-u}(F))^{-1})^{1/2}.
\end{equation} 
Similarly, we can define  \begin{equation} \label{eq:rhoe}
\rho_{2k}^{or}(x)=(\prod_{u^{2k}-1=0 \atop u\neq 1}
\gamma_{u/u-1}(x))^{1/2}. \end{equation}

L{\scriptsize{EMMA}} 2.2
({\it An analogue of Proposition~$5.3$ of} \cite{ati-tal}).
{\it If $(p,k)=1$  then the following 
diagram where $i(x)=1\otimes x$ is commutative:}

\begin{center}
\begin{picture}(350,88)
\put(50,72){$\widetilde{KSO}(X)$}
\put(70,67){\vector(0,-1){40}}
\put(50,12){$\widetilde{KSO}(X)_p$}
\put(75,47){$i$}
\put(95,67){\vector(3,-2){60}}
\put(130,47){$\theta_{k}^{or}$}
\put(145,12){$1+\widetilde{KSO}(X)_p$.}
\put(100,17){\vector(1,0){43}}
\put(115,22){$\rho_{k}^{or}$}
 \end{picture} \end{center}

R{\scriptsize EMARK}. If $(p,k)=1,$ then  ${\Bbb Q}_k \subseteq {\Bbb Z}_p.$ So,
 using Proposition~3.15.2 of \cite{die} and Examples~ 5.14 and 5.15
of \cite{ada}-II, we see that $\theta_{k}^{or}$ agrees with
Bott operation $\theta_k$ which is denoted by $\rho^k$ in \cite{ada}-II. 

Now, we give our main theorem.

T{\scriptsize HEOREM} 2.3. {\it
Let $p$ be a prime number and $k_p$ be a generator of
 $({\Bbb Z}/p^2{\Bbb Z})^*,$ the group of units in
${\Bbb Z}/p^2{\Bbb Z}$. Then the following 
diagram is commutative}:

\begin{picture}(358,118)
\put(50,105){$0$}
\put(53,102){\vector(0,-1){25}}
\put(35,63){$\widetilde{KSO}(X)_{(p),\Gamma }$}
\put(105,65){\vector(1,0){55}}
\put(130,73){$\tilde{q}$}
\put(165,63){$\widetilde{KSO}(X)_{(p)}/TO(X)_{(p)}$}
\put(285,65){\vector(1,0){40}}
\put(330,61){$0$}
\put(218,55){\vector(0,-1){35}}
\put(222,40){$\tilde{\theta}_{k_p}^{or}$}
\put(53,60){\vector(0,-1){40}}
\put(60,40){$i_{\Gamma}$}
\put(3,5){$0$}
\put(10,9){\vector(1,0){25}}
\put(35,5){$\widetilde{KSO}(X)_{p,\Gamma }$}
\put(100,7){\vector(1,0){75}}
\put(130,15){$\rho_{k_p,\Gamma}^{or}$}
\put(180,5){$1+\widetilde{KSO}(X)_{p,\Gamma }$}
\put(265,7){\vector(1,0){60}}
\put(330,3){$0$}
\end{picture}\\
{\it Here the index $\Gamma$ indicates that we factor out the image of
$(\psi^{k_p}-1)$ and $\tilde{q}$ is the quotient map}.

P{\scriptsize ROOF}. First, we show that rows and columns are well-defined and
exact.\\
(a) Using Lemma~2.1, the fact that localization and completion are exact
functors on the category of finitely generated abelian groups and the 
naturality of Adams'~ operations, we have the following identifications:
$$\widetilde{KSO}(X)_{(p)}/(\psi^{k_p}-1)(\widetilde{KSO}(X)_{(p)})=
\widetilde{KSO}(X)_{(p)}/((\psi^{k_p}-1)(\widetilde{KSO}(X)))_{(p)}= $$
$$(\widetilde{KSO}(X)/(\psi^{k_p}-1)(\widetilde{KSO}(X)))_{(p)} \subseteq
\widetilde{KSO}(X)_p/((\psi^{k_p}-1)(\widetilde{KSO}(X)))_p= $$
$$\widetilde{KSO}(X)_p/(\psi^{k_p}-1)(\widetilde{KSO}(X)_p).$$
Hence, $i:\widetilde{KSO}(X)_{(p)}\rightarrow \widetilde{KSO}(X)_p $
defined by  $i(x/m)=(m)^{-1}\otimes x$ induces a monomorphism 
$i_{\Gamma}:\widetilde{KSO}(X)_{(p),\Gamma}
\rightarrow \widetilde{KSO}(X)_{p,\Gamma}.$\\
(b) By Theorem~4.5 of Atiyah-Tall \cite{ati-tal}, $\rho_{k_p}^{or}$
 induces an isomorphism
$$\rho_{k_p,\Gamma}^{or}:\widetilde{KSO}(X)_{p,\Gamma}\rightarrow
 1+\widetilde{KSO}(X)_{p,\Gamma}.$$ 
(c) To show that $(\psi^{k_p}-1)(KSO(X)_{(p)}) \subseteq TO(X)_{(p)}.$ Let
 $$\frac{E-F}{m} \in \widetilde{KSO}(X)_{(p)}.$$ Then 
$$ (\psi^{k_p}-1)(\frac{E-F}{m})=(\frac{\psi^{k_p}E-E}{m})-
(\frac{\psi^{k_p}F-F}{m}).$$
By Quillen \cite{qui}, there is a fiberwise
 map of degree a power of $k_p$ between $\psi^{k_p}E$ and $E.$ So,
 by Dold's Theorem
mod $k$ \cite{ada}-I we have $$k_{p}^{e}(\psi^{k_p}E-E) \in TO(X)$$
 for some integer $e.$ Since $(p,k_p)=1,$ we have
$$(\frac{\psi^{k_p}E-E}{m})=\frac{k_{p}^{e}(\psi^{k_p}E-E)}
{k_{p}^{e} m} \in TO(X)_{(p)}.$$
Similarly, $$\frac{\psi^{k_p}F-F}{m} \in TO(X)_{(p)} $$ and hence
$$ (\psi^{k_p}-1)(\frac{E-F}{m})\in TO(X)_{(p)}.$$ Thus,
 we have an epimorphism
$\tilde{q}:\widetilde{KSO}(X)_{(p),\Gamma}\rightarrow 
\widetilde{KSO}(X)_{(p)}/TO(X)_{(p)}.$ \\
(d) It is easy to see that
$\theta_{k_p}^{or}:\widetilde{KSO}(X)_{(p)}\rightarrow 1+\widetilde{KSO}(X)_p$ 
given by $\theta_{k_p}^{or}(x/m)=(\theta_{k_p}^{or}(x))^{1/m}$
 is an exponential
map. Let $$\frac{E-F}{m} \in TO(X)_{(p)}.$$ Then $nS(E)$ is stably 
fibre homotopy equivalent to $nS(F)$ for some $n$ with $(p,n)=1.$ 
So by \cite{ada}-(II)  Corollary~5.8, $$ \theta_{k_p}^{or}(E-F)^n=
\frac{\psi^{k_p}(1+u)}{1+u} \mbox{ in }1+\widetilde{KSO}(X)_p $$ for some 
$u\in \widetilde{KSO}(X).$ Since $(p,n)=1,$  $(1+u)^{1/nm}=1+w$ in
$1+\widetilde{KSO}(X)_p$ for some $w\in \widetilde{KSO}(X)_p.$ Hence
$$\theta_{k_p}^{or}(\frac{E-F}{m}) =\theta_{k_p}^{or}(E-F)^{1/m}=
(\theta_{k_p}^{or}(E-F)^n)^{1/nm}$$ 
$$=\frac{\psi^{k_p}(1+u)^{1/nm}}{(1+u)^{1/nm}}=
\frac{\psi^{k_p}(1+w)}{1+w}.$$ Thus $\theta_{k_p}^{or}$
 induces a homomorphism
$$\tilde{\theta}_{k_p}^{or}:\widetilde{KSO}(X)_{(p)}/TO(X)_{(p)}\rightarrow  
1+\widetilde{KSO}(X)_{p,\Gamma}.$$
Finally, we show the commutativity of our diagram.\\
Let $x/m \in \widetilde{KSO}(X)_{(p)}.$ Then  $\tilde{\theta}_{k_p}^{or} \circ
\tilde{q}(x/m+(\psi^{k_p}-1)(\widetilde{KSO}(X)_{(p)}))=$
 $\tilde{\theta}_{k_p}^{or}(x/m+TO(X)_{(p)})=$
$\theta_{k_p}^{or}(x)^{1/m}+(\psi^{k_p}-1)(1+\widetilde{KSO}(X)_p).$ On the 
other hand, 
$(\rho_{k_p,\Gamma}^{or} \circ i_{\Gamma})(x/m+(\psi^{k_p}-1)
(\widetilde{KSO}(X)_{(p)}))=$ 
$\rho_{k_p,\Gamma}^{or} (i(x/m)+(\psi^{k_p}-1)(\widetilde{KSO}(X)_p))=$
$\rho_{k_p}^{or}(i(x/m))+(\psi^{k_p}-1)(1+\widetilde{KSO}(X)_p)).$ Now, the
 result follows from Lemma~2.2. This completes the proof of Theorem~{2.3}.

C{\scriptsize OROLLARY} 2.4 ({\it Formula~I of $TO(X)_{(p)}$}).
$$ TO(X)_{(p)}=(\psi^{k_p}-1)(\widetilde{KSO}(X)_{(p)}).$$

P{\scriptsize ROOF}. Since  $\tilde{\theta}_{k_p}^{or} \circ \tilde{q}=
 \rho_{k_p,\Gamma}^{or} \circ i_{\Gamma},$  $\tilde{q}$ is injective and
hence an isomorphism. So, 
$TO(X)_{(p)}=(\psi^{k_p}-1)(\widetilde{KSO}(X)_{(p)}).$

C{\scriptsize OROLLARY} 2.5 ({\it Formula~II of $TO(X)_{(p)}$}).
\begin{eqnarray}
TO(X)_{(p)}&=\{& x \in 
\widetilde{KSO}(X)_{(p)}:\theta_{k_p}^{or}(x)=\frac{\psi^{k_p}(1+u)}{1+u} 
\mbox{ in } 1+\widetilde{KSO}(X)_{p} \nonumber \\ & & \mbox{ for some } u 
\in \widetilde{KSO}(X)_{p}\}. \nonumber \end{eqnarray}

P{\scriptsize ROOF}. Clearly, the right hand side of the above equality
 is a well-defined subgroup of $\widetilde{KSO}(X)_{(p)}.$ The fact that
 $i_{\Gamma}$ is injective implies that
\begin{equation} \label{eqn:J1}
 i( \widetilde{KSO}(X)_{(p)})\cap (\psi^{k_p}-1)
(\widetilde{KSO}(X)_p)=i((\psi^{k_p}-1)(\widetilde{KSO}(X)_{(p)})). 
\end{equation}
The fact that $\rho_{k_p,\Gamma}^{or}$ is an isomorphism implies that 
\begin{equation} \label{eqn:J2}
\rho_{k_p}^{or} (\psi^{k_p}-1)(\widetilde{KSO}(X)_p)=
(\psi^{k_p}-1)(1+\widetilde{KSO}(X)_p). \end{equation}
 Now let $x \in TO(X)_{(p)} $ then by  Formula~I of
 $TO(X)_{(p)},$ 
$x \in (\psi^{k_p}-1)\widetilde{KSO}(X)_{(p)}.$ Hence 
from (\ref{eqn:J1}) and (\ref{eqn:J2})
$$\theta_{k_p}^{or}(x)=\rho_{k_p}^{or}(i(x))=\frac{\psi^{k_p}(1+u)}{1+u}
\mbox{ in } 1+\widetilde{KSO}(X)_p $$ for some $ u
\in \widetilde{KSO}(X)_p.$

If $X$ is a finite $CW$ complex, then $\widetilde{JO}(X)$ is a finite abelian
group. So, by Lemma~2.1, the $p$-primary factor of the order of 
$x+TO(X)\in \widetilde{JO}(X)$ is the order of  
$x+TO(X)_{(p)}\in \widetilde{JO}(X)_{(p)},$  the smallest
$p^m$ such that $p^mx\in TO(X)_{(p)}.$ 
\vspace{1cm}\\
{\large \textbf{3. $J$-orders of elements of $\widetilde{KO}({\Bbb C}  P^m)$ }}

We will show how to use Formulae I and II of $TO(X)_{(p)}$ to
find the $J$-orders of elements of $\widetilde{KO}({\Bbb C}  P^m).$ As we have
shown in \cite{obi}, we only need to consider the case $m$ is even, that is
 $m=2t$ for some $t\in {\Bbb N}.$ Let  
 $ P_{m}(y;m_1,\ldots,m_t)=m_1y+m_2 y^2+\cdots+m_ty^t \in
\widetilde{KO}({\Bbb C}P^m)={\Bbb Z}[y]\mbox{ (mod } y^{t+1}).$ In order to
find  $b_m(P_{m}(y;m_1,\ldots,m_t)),$ the $J$-order
of $ P_{m}(y;m_1,\ldots,m_t),$ the following two lemmas will be useful.

L{\scriptsize EMMA} 3.1. 
{\it Let $k_p$ be a generator of $({\Bbb Z}/p^2{\Bbb Z})^*$.
 If n $\in {\Bbb N},$ then}
\begin{enumerate}
\item[\rm(i)] $\nu_2(3^{2n}-1)=3+\nu_2(n).$ 
\item[\rm(ii)] {\it For an odd prime $p,$}
$$\nu_p(k_{p}^{2n}-1)=\left\{ \begin{array}{ccc}
0&\mbox{if}& 2n\not\equiv 0 \mbox{ mod}(p-1)  \\
1+\nu_p(n) & \mbox{if}& 2n\equiv 0 \mbox{ mod}(p-1). \end{array}
\right.$$
\end{enumerate}

P{\scriptsize ROOF}. (i) is well-known.\\
(ii) Let  $\nu_p(k_{p}^{2n}-1)=s.$ Then $k_{p}^{2n}\equiv 1 \mbox{ mod }
p^s.$ If $s\geq 1,$ then $({\Bbb Z}/p^s{\Bbb Z})^*$ is cyclic of order
$p^{s-1}(p-1)$ with generator $k_p$ (\cite{ire-ros}, Theorem~2, p. 43).
 So, $2n=p^{s-1}(p-1)d$ for some
$d\in {\Bbb N}$ with $(d,p)=1$(\cite{ire-ros}, Lemma~3, p. 42).
 Hence, $s=1+\nu_p(n).$

L{\scriptsize EMMA} 3.2.
{\it Let $k_p$ be a generator of $({\Bbb Z}/p^2{\Bbb Z})^*$ and 
$r,s \in {\Bbb N}$ with $r\geq s.$ Then}
\begin{enumerate}
\item[\rm(i)] $\nu_2(\prod_{i=s}^{r}(3^{2i}-1))=3(r-s+1)+
\sum_{i=s}^{r}\nu_2(i).$  
\item[\rm(ii)] {\it For an odd prime} $p,\,
\nu_p(\prod_{i=s}^{r}(k_{p}^{2i}-1))$ $$= \left\{ \begin{array}{ccc}
[\frac{2r}{p-1}]+\sum_{i=1}^ {[\frac{2r}{p-1}]}\nu_p(i)-
[\frac{2(s-1)}{p-1}]-\sum_{i=1}^ {[\frac{2(s-1)}{p-1}]}\nu_p(i)& \mbox{if}&
p\leq2r+1\\
0&\mbox{if}& p>2r+1. \end{array}\right.$$
\end{enumerate}

P{\scriptsize ROOF}. (i) $\nu_2(\prod_{i=s}^{r}(3^{2i}-1))=
\sum_{i=s}^{r}\nu_2(3^{2i}-1)= \sum_{i=s}^{r}(3+\nu_2(i))=
3(r-s+1)+\sum_{i=s}^{r}\nu_2(i).$\\
(ii) If $p>2r+1,$ then $\nu_p(k_{p}^{2i}-1)=0$ for each $i=s,\ldots,r.$
 Hence
$$\nu_p(\prod_{i=s}^{r}(k_{p}^{2i}-1))=0.$$ 
If $p\leq 2r+1,$ then $p-1=2d$ for some $d\in \{1,\ldots,r\}.$
$$\nu_p(\prod_{i=s}^{r}(k_{p}^{2i}-1))=\sum_{i=s}^{r}\nu_p(k_{p}^{2i}-1)=$$
$$ \sum_{i=s \atop 2i\equiv 0 \mbox{ {\tiny mod}}(p-1)}^ {r}(1+\nu_p(i))=
\sum_{i=1}^ {[\frac{2r}{p-1}]}(1+\nu_p(2di))-
\sum_{i=1}^ {[\frac{2(s-1)}{p-1}]}(1+\nu_p(2di)).$$ But $\nu_p(2di)=
\nu_p(i).$ So $$\nu_p(\prod_{i=s}^{r}(k_{p}^{2i}-1))=
\sum_{i=1}^ {[\frac{2r}{p-1}]}(1+\nu_p(i))
-\sum_{i=1}^ {[\frac{2(s-1)}{p-1}]}(1+\nu_p(i))=$$
$$[\frac{2r}{p-1}]+\sum_{i=1}^ {[\frac{2r}{p-1}]}\nu_p(i)-
[\frac{2(s-1)}{p-1}]-\sum_{i=1}^ {[\frac{2(s-1)}{p-1}]}\nu_p(i).$$
This completes the proof.

Now, let $k_p$ be an odd generator of $({\Bbb Z}/p^2{\Bbb Z})^*,$ say
$k_p=2q+1$ (take $k_2=3$).

R{\scriptsize EMARK}. We take $k_p$ to be odd only to reduce the work, even
 $k_p$ 
 works equally well.

According to Formula~I of $TO(X)_{(p)},$
 $\nu_p(b_m(P_{m}(y;m_1,\ldots,m_t)))$ is the smallest non-negative 
integer $v$ such that
 \begin{equation} \label{eq:f1}
p^vP_{m}(y;m_1,\ldots,m_t)=(1-\psi^{k_p})(u) 
\end{equation}
in $\widetilde{KO}({\Bbb C}  P^m)_{(p)}$ for some $u \in 
\widetilde{KO}({\Bbb C}  P^m)_{(p)}.$

From \cite{ada-wal} Theorem~2.2, and \cite{mps} Lemma~3.6,
 \begin{equation} \label{eq:h1} 
\psi^{k_p}(y)=y(\sum_{j=0}^{q}
\frac{k_p}{2j+1}{q+j \choose 2j}y^j)^2=y(\sum_{j=0}^{q}b_jy^j)^2
\end{equation}
where $$b_j=\frac{k_p}{2j+1}{q+j \choose 2j},\,j=0,\ldots,q.$$ So, for
$r=2,\ldots,t$
$$\psi^{k_p}(y^r)=(\psi^{k_p}(y))^r=\sum_{j=0 \atop j\leq t-r}^{2rq}
C_{j,r}y^{r+j}$$ where
$$ C_{j,r}=\sum_{i_1+\cdots+i_{2r}=j \atop i_1,\ldots,i_{2r}\in 
\{0,\ldots,q\}}^{}b_{i_1}b_{i_2}\ldots b_{i_{2r}}.$$

Let $u\in\widetilde{KO}({\Bbb C}  P^m)_{(p)}.$ Then $u=a_1y+\cdots+
a_ty^t$ for some $a_i\in {\Bbb Z}_{(p)}.$ Using (\ref{eq:h1}), it is
easy to see that the coefficient of $y^r$ in $(1-\psi^{k_p})(u)$ is
$$-\sum_{i=j_r}^{r-1}C_{r-i,i}a_i+(1-k_{p}^{2r})a_r$$ where
$j_r=[\frac{r-1}{k_p}]+1.$\\
So, from (\ref{eq:f1}), we need to find the smallest $v$ which solves
the following system of equations in ${\Bbb Z}_{(p)}:$
$$-\sum_{i=j_r}^{r-1}C_{r-i,i}a_i+(1-k_{p}^{2r})a_r=p^vm_r $$
where $r=1,\ldots,t.$\\
The above system has the following solutions:
$$a_r=\frac{p^vM_{k_p,r}(m_1,\ldots,m_t)}{(1-k_{p}^{2})
\ldots (1-k_{p}^{2r})}$$ where $M_{k_p,1}(m_1,\ldots,m_t)=m_1 $ and
for $r=2,\ldots,t$
$$M_{k_p,r}(m_1,\ldots,m_t)=\sum_{i=j_r}^{r-1}C_{r-i,i}(1-k_{p}^{2(i+1)})
\ldots (1-k_{p}^{2(r-1)})M_{k_p,i}(m_1,\ldots,m_t)$$ 
$$+m_r(1-k_{p}^{2})\ldots (1-k_{p}^{2(r-1)}).$$ 
Now, $a_r\in {\Bbb Z}_{(p)}$ implies that $\nu_p(a_r)\geq 0.$ So
$$ v\geq\,\max_{ r=1,\ldots,t}\,\{\,\nu_p\,(\prod_{i=1}^{r}(1\!-\!
k_{p}^{2i}))\!-\!\nu_p(M_{k_p,r}(m_1,\ldots,m_t)),0:
\! M_{k_p,r}(m_1\!,\!\ldots\!,\!m_t)\neq 0\}$$ 
is a necessary and sufficient condition on $v$ so that (\ref{eq:f1})
is satisfied. Hence, we have:

T{\scriptsize HEOREM} 3.3.
 $\nu_p(b_m(P_{m}(y;m_1,\ldots,m_t)))=$
$$\max_{ r=1,\ldots,t}\,\{\,\nu_p\,(\prod_{i=1}^{r}(1\!-\!
k_{p}^{2i}))\!-\!\nu_p(M_{k_p,r}(m_1,\ldots,m_t)),0:
\! M_{k_p,r}(m_1\!,\!\ldots\!,\!m_t)\neq 0\}.$$ 

Now, let us use Formula~II.\\
Let $\theta_{k_p}(P_{m}(y;m_1,\ldots,m_t))=1+\alpha_1(m_1,\ldots,m_t)y+
\cdots+\alpha_t(m_1,\ldots,m_t)y^t$ for some 
$\alpha_i(m_1,\ldots,m_t)\in {\Bbb Z}_p$
 (see \cite{obi}, Theorem~2.2). $\nu_p(b_m(P_{m}(y;m_1,\ldots,m_t)))$
 is the smallest non-negative integer $v$ such that
 \begin{equation} \label{eq:f3}
\theta_{k_p}(P_{m}(y;m_1,\ldots,m_t))^{p^v}=\frac{\psi^{k_p}(1+u)}{1+u}
\mbox{ in } 1+\widetilde{KO}({\Bbb C}  P^m)_{p}\end{equation}
for some $u\in \widetilde{KO}({\Bbb C}  P^m)_{p}.$
Let $u=b_1y+\cdots+b_ty^t$ for some $b_i\in {\Bbb Z}_p.$
With the above symbols, the coefficient of $y^r$ in $\psi^{k_p}(u)$ is
$$\sum_{i=j_r}^{r-1}C_{r-i,i}b_i+b_rk_{p}^{2r}.$$ To avoid excessive notation,
we write
$\theta_{k_p}(P_{m}(y;m_1,\ldots,m_t))^{p^v}=1+\alpha_1y+\cdots+\alpha_ty^t$
 where $\alpha_i$  inolves quantities containing $p$ in some way.

From (\ref{eq:f3}), we have  $1+\psi^{k_p}(u)=1+d_1y+\cdots+d_ty^t$ where
$$d_n=\sum_{{i+s=n \atop b_0=\alpha_0=1}}b_i\alpha_s.$$ Thus
$$ \sum_{i=j_r}^{r-1}C_{r-i,i}b_i+b_rk_{p}^{2r}=b_r+\sum_{{s>0 \atop i+s=r}}
b_i\alpha_s$$ which implies that
$$b_r=\frac{\sum_{{s>0 \atop i+s=r}}b_i
\alpha_s-\sum_{i=j}^{r-1}C_{r-i,i}b_i}{(k_{p}^{2r}-1)}=
\frac{L_{k_p,r}(m_1,\ldots,m_t)}{(k_{p}^{2}-1)\ldots(k_{p}^{2r}-1)}$$ 
where $L_{k_p,1}(m_1,\ldots,m_t)=\alpha_1$ and for $r=2,\ldots,t,$
 $$L_{k_p,r}(m_1,\ldots,m_t)=\sum_{i=1}^{r-1}L_{k_p,r-i}(m_1,\ldots,m_t)
\alpha_i(k_{p}^{2(r-i+1)}-1)\ldots
(k_{p}^{2(r-1)}-1)-$$ $$\sum_{i=j_r}^{r-1}C_{r-i,i}L_{k_p,i}(m_1,\ldots,m_t)
(k_{p}^{2(i+1)}-1)
\ldots(k_{p}^{2(r-1)}-1)+\alpha_r(k_{p}^{2}-1)\ldots(k_{p}^{2(r-1)}-1).$$
Now, $b_i\in{\Bbb Z}_p\mbox{ for } i=1,\ldots,t$ implies that
$$\nu_p(L_{k_p,r}(m_1,\ldots,m_t))\geq \nu_p((k_{p}^{2}-1)\ldots
(k_{p}^{2r}-1)).$$ So, we have:

T{\scriptsize HEOREM} 3.4.
$\nu_p(b_m(P_{m}(y;m_1,\ldots,m_t)))$ {\it is the smallest $v$ such that}\\
$\nu_p(L_{k_p,r}(m_1,\ldots,m_t))\geq \nu_p((k_{p}^{2}-1)\ldots
(k_{p}^{2r}-1))$ {\it  for  each } $r=1,\ldots,t.$

 Using Lemma~3.2 and any one of the above two theorems, we directly obtain:

C{\scriptsize OROLLARY} 3.5.
{\it If $p>2t+1,$ then}
$\nu_p(b_m(P_{m}(y;m_1,\ldots,m_t)))=0.$ {\it Consequently},
$$\widetilde{JO}({\Bbb C}  P^m)\cong \bigoplus_{\mbox{for all
 primes } p\leq m+1 }\widetilde{JO}({\Bbb C}  P^m)_{(p)}.$$

From Theorem~3.3, to find $b_m(P_{m}(y;m_1,\ldots,m_t))$ we only need to
find \\ $\nu_p(M_{k_p,r}(m_1,\ldots,m_t))$ for $r=1,\ldots,t.$ Therefore,
 it may be a good problem if one tries to obtain a general formula for 
$\nu_p(\!M_{k_p,r}(m_1,\ldots,m_t)\!)$ in term of $r,k_p,m_1,\ldots,m_t.$\\
Next, we compute $\nu_p(M_{k_p,r}(0,\ldots,m_r=1,0,\ldots,0))$ for
$p=2,3$  and then we obtain  simple formulae for the $2$ and $3$ primary
factors of the $J$-orders of the canonical generators of
$\widetilde{JO}({\Bbb C}  P^m).$ These simple formulae have been already
conjectured in \cite{obi}.

For $n=1,\ldots,t,$ the $J$-order of $y^n+TO({\Bbb C}P^m)$ is 
$b_m(P_m(y;0,\ldots,m_n=1,0,\ldots,0)).$ Let $M_{k_p,r}=M_{k_p,r}(0,\ldots,
m_n=1,0,\ldots,0)/(1-k_{p}^{2})\cdots(1-k_{p}^{2(n-1)}).$ 
Then $M_{k_p,r}=0$ for $r<n,$ $M_{k_p,n}=1$  and for $r=n+1,\ldots,t,$
$$M_{k_p,r}=\sum_{{i=j_r \atop i\geq n}}^{r-1}C_{r-i,i}(1-k_{p}^{2(i+1)})
\ldots(1-k_{p}^{2(r-1)})M_{k_p,i}$$ where
$j_r=[\frac{r-1}{k_p}]+1$\\
 Hence, from Theorem~3.3, we have
$$\nu_p(b_m(y^n))=\max_{r=n,\ldots,t}\{\nu_p(\prod_{i=n}^{r}(1-k_{p}^{2i}))-
\nu_p(M_{k_p,r}),0: M_{k_p,r}\neq 0\}.$$

P{\scriptsize ROPOSITION} 3.6.
{\it If $p=2\mbox{ or }3,$ then} 
$$\nu_p(M_{k_p,r})=\sum_{s=1}^{[\frac{2(r-1)}{p-1}]}
\nu_p(s)-\sum_{s=1}^{[\frac{2(n-1)}{p-1}]}\nu_p(s)$$
{\it for each} $r=n,\ldots,t.$

P{\scriptsize ROOF}. We prove this proposition for $p=2$ (the case $p=3$ is 
similar). Recall that $k_2=3.$ So we need to show that
$\nu_2(M_{3,r})=r-n+\sum_{s=n}^{r-1}\nu_2(s)$ for $r=n+1,\ldots,t,$ where
$$M_{3,r}=\sum_{i=j_r \atop i\geq n}^{r-1}3^{3i-r}{2i\choose r-i}(1-3^{2(i+1)})
\cdots (1-3^{2(r-1)})M_{3,i}$$
 by induction on $r.$ If 
$r=n+1$ then $\nu_2(M_{3,r})=\nu_2{2n \choose 1}=1+\nu_2(n)$.
 So let $n+1<r\leq t.$ We claim that
$\nu_2(3^{3i-r}{2i \choose r-i}(1-3^{2(i+1)})\cdots (1-3^{2(r-1)})
M_{3,i})>\nu_2(3^{2r-3}{2(r-1) \choose 1}M_{3,r-1})$ for each
$\max\{j_r,n\}\leq i <r-1.$
 Suppose that $\max\{j_r,n\}\leq i<r-1.$  Then
 by induction hypothesis and Lemma~3.2, 
$$\nu_2(3^{3i-r}{2i \choose r-i}(1-3^{2(i+1)})\cdots (1-3^{2(r-1)})
M_{3,i})$$ $$=\nu_2{2i \choose r-i}+3(r-i-1)+i-n-\nu_2(i)+
\sum_{s=n}^{r-1}\nu_2(s).$$ On the other hand,
$$\nu_2(3^{2r-3}{2(r-1) \choose 1}M_{3,r-1})=r-n+\sum_{s=n}^{r-1}\nu_2(s). $$
So, we need to show that
 $\nu_2{2i \choose r-i}+2(r-i-1)>\nu_2(2i).$ But this follows directly from
the fact that $\nu_2{2i\choose r-i}=\nu_2(2i)-\nu_2(r-i)$ if
$\nu_2(2i)\geq r-i-1.$ This completes the proof of our claim. Hence,
 $\nu_2(M_{3,r})=\nu_2(3^{2r-3}{2(r-1) 
\choose 1}M_{3,r-1})=r-n+\sum_{s=n}^{r-1}\nu_2(s).$ This completes the proof.

Unfortunately, the above proof can not be used for $p\neq 2,3.$

T{\scriptsize HEOREM} 3.7.
{\it If $p=2\mbox{ or }3$ and $1\leq n \leq t.$ Then}
$$\nu_p(b_m(y^n))=\max \left\{ s-[\frac{2(n-1)}{p-1}]+\nu_p(s):
[\frac{2n}{p-1}]\leq s\leq[\frac{2t}{p-1}] \right\}. $$

P{\scriptsize ROOF}. Let $p=2,$ then
 $\nu_2(b_m(y^n))=\max\{\nu_2(\prod_{i=n}^{r}(1-3^{2i}))-
\nu_2(M_{3,r}):r=n,\ldots,t\}=$
$\max\{2r-2n+2+\nu_2(2r):r=n,\ldots,t\}=$  
$\max\{ s-2(n-1)+\nu_2(s):2n\leq s\leq 2t \}.$
The case $p=3$ is similar.

R{\scriptsize EMARK}. If Proposition~3.6 is true for 
some values of $p$ other
than $2$ or $3,$ then Theorem~3.7 will be also true for those values of $p.$
\vspace{1cm}\\
{\large \textbf{4. An illustrative example $\widetilde{JO}({\Bbb C}P^4)$}}

If $\widetilde{KO}(X)\!=\!<y_1,\ldots,y_n>,$ then $\widetilde{JO}(X)_{(p)}=
<\!\alpha_{1,p}\!=y_1+TO(X)_{(p)},\ldots,\alpha_{n,p}\!=y_n+TO(X)_{(p)}>.$
 So to compute
$\widetilde{JO}(X)_{(p)},$ we need to find all relations between
$\alpha_{1,p},\ldots,\alpha_{n,p},$ i.e., we need to find 
``sufficient'' solutions for the equation:
\begin{equation} \label{eq:f4} c_1\alpha_{1,p}+\cdots+c_n\alpha_{n,p}=0
 \mbox{  in  } 
\widetilde{JO}(X)_{(p)},\, c_1,\ldots,c_n \in {\Bbb Z}. \end{equation}
This 
implies that $c_1y_1+\cdots+c_ny_n \in TO(X)_{(p)}.$
 Now using formulae I and II
of $TO(X)_{(p)},$ one may try to find ``sufficient'' solutions for
(\ref{eq:f4}).

$\widetilde{KO}({\Bbb C}  P^4)=\{a_1y+a_2y^2:a_1,a_2 \in {\Bbb Z},\,
 y^3=0\}.$ So, $\widetilde{JO}({\Bbb C}  P^4)_{(p)}=
<\alpha_{1,p}=y+TO({\Bbb C}P^4)_{(p)},\,\alpha_{2,p}=
y^2+TO({\Bbb C}P^4)_{(p)}>
=<\alpha_{1,p}>+<\alpha_{2,p}>.$ To find relations between $\alpha_{1,p}$ and 
$\alpha_{2,p},$ we need to solve $c_1\alpha_{1,p}+c_2\alpha_{2,p}=0$ in
$\widetilde{JO}({\Bbb C}  P^4)_{(p)}.$

 $\widetilde{JO}({\Bbb C}  P^4)_{(2)}=
<\alpha_{1,2}=y+TO({\Bbb C}P^4)_{(2)}>+<\alpha_{2,2}=
y^2+TO({\Bbb C}P^4)_{(2)}>.$ 
 $<\alpha_{1,2}>$ is cyclic of order $64$ and $<\alpha_{2,2}>$ 
is cyclic of order $16.$ Also, $2\alpha_{2,2}=40\alpha_{1,2}.$ Hence
$\widetilde{JO}({\Bbb C}P^4)_{(2)}\cong {\Bbb Z}/2{\Bbb Z}\oplus
{\Bbb Z}/64{\Bbb Z}.$ Similarly, 
$\widetilde{JO}({\Bbb C}  P^4)_{(3)}\cong {\Bbb Z}/9{\Bbb Z}$ and 
$\widetilde{JO}({\Bbb C}  P^4)_{(5)}\cong {\Bbb Z}/5{\Bbb Z}.$
 Thus, by Corollary~3.5, we get a known result:

T{\scriptsize HEOREM} 4.1.
$\widetilde{JO}({\Bbb C}  P^4)\cong {\Bbb Z}/2{\Bbb Z}\oplus
{\Bbb Z}/64{\Bbb Z}\oplus
{\Bbb Z}/9{\Bbb Z}\oplus {\Bbb Z}/5{\Bbb Z}.$

$\begin{array}{c}\mbox{ Middle East Technical University } \\  
\mbox{ Department of Mathematics }\\
\mbox{ Ankara 06531, TURKEY } \\
\mbox{ e-mail : b092385@rorqual.cc.metu.edu.tr  } \end{array}$
 \end{document}